\documentclass[12pt]{article}
\usepackage{amsmath}
\usepackage{amssymb}
\usepackage[dvips]{graphicx}
\usepackage{epsfig}
\def\eps{\varepsilon}
\font\tencmmib=cmmib10 \skewchar\tencmmib '60
\newfam\cmmibfam
\textfont\cmmibfam=\tencmmib

\def\bbox{\quad\hbox{\vrule \vbox{\hrule \vskip2pt \hbox{\hskip2pt
\vbox{\hsize=1pt}\hskip2pt} \vskip2pt\hrule}\vrule}}
\def\lessim{\ \lower4pt\hbox{$
\buildrel{\displaystyle <}\over\sim$}\ }
\def\gessim{\ \lower4pt\hbox{$\buildrel{\displaystyle >}
\over\sim$}\ }

\def\eps{{\varepsilon}}

\def\la{{\Bigl\langle}}
\def\ra{{\Bigr\rangle}}

\def\qed{\hfill\break\rightline{$\bbox$}}
\newcommand{\E}{\mathbb{E}}

\newcommand{\Reals}{\mathbb{R}}

\newtheorem{theorem}{Theorem}

\makeatletter
\@addtoreset{equation}{section}

\makeatother

\font\tencmmib=cmmib10 \skewchar\tencmmib '60
\newfam\cmmibfam
\textfont\cmmibfam=\tencmmib

\def\bbox{\quad\hbox{\vrule \vbox{\hrule \vskip2pt \hbox{\hskip2pt
\vbox{\hsize=1pt}\hskip2pt} \vskip2pt\hrule}\vrule}}
\def\lessim{\ \lower4pt\hbox{$
\buildrel{\displaystyle <}\over\sim$}\ }
\def\gessim{\ \lower4pt\hbox{$\buildrel{\displaystyle >}
\over\sim$}\ }

\def\eps{\varepsilon}

\def\go0{\to 0}

\def\la{\langle}

\def\leftitem#1{\item{\hbox to\parindent{\enspace#1\hfill}}}

\def\qed{\hfill\break\rightline{$\bbox$}}

\def\ra{\rangle}

\def\sg{\sigma}

\def\sg2{\sigma^2}

\def\__{_{\infty}}

\begin{document}

\title{A deletion-invariance property for random measures satisfying the Ghirlanda-Guerra identities.}
\author{Dmitry Panchenko\thanks{Department of Mathematics, Texas A\&M University, email: panchenk@math.tamu.edu. Partially supported by NSF grant.}}

\maketitle
\begin{abstract} 
We show that if a discrete random measure on the unit ball of a separable Hilbert space 
satisfies the Ghirlanda-Guerra identities then by randomly deleting half of the points 
and renormalizing the weights of the remaining points we obtain the same random
measure in distribution up to rotations.
\end{abstract}
\vspace{0.5cm}
Key words: Sherrington-Kirkpatrick model, Gibbs measure.

\section{Introduction and main result.} 
Let us consider a countable index set $A$ and random probability measure $\mu$ 
on a unit ball $B$ of a separable Hilbert space $H$ such that 
$\mu=\sum_{\alpha\in  A}w_\alpha \delta_{\xi_\alpha}$ 
for some random points $\xi_\alpha \in B$ and weights $(w_\alpha)$.
We will call indices $\alpha$ from the set $A$ ``configuration'' and
for a function $f=f(\alpha^1,\ldots,\alpha^n)$ of $n$ configurations we will denote its average
with respect to the measure $\mu$ by
\begin{equation}
\la f \ra = \sum_{\alpha^1,\ldots, \alpha^n} w_{\alpha^1}\cdots w_{\alpha^n} f(\alpha^1,\ldots,\alpha^n).
\label{Av}
\end{equation}
We say that the random measure $\mu$ satisfies the Ghirlanda-Guerra identities \cite{GG}
if for any $n\geq 2$ and any function $f$ that depends on the configurations $\alpha^1,\ldots,\alpha^n$ 
only through the scalar products, or overlaps, $R_{\ell,\ell'} = \xi_{\alpha^\ell}\cdot \xi_{\alpha^{\ell'}}$ 
for $\ell,\ell^\prime\leq n$ we have
\begin{equation}
\E \la f R_{1,n+1}^p \ra = \frac{1}{n}\, \E \la f \ra \, \E \la R_{1,2}^p \ra + 
\frac{1}{n} \sum_{\ell=2}^{n} \E \la f R_{1,\ell}^p\ra
\label{GG}
\end{equation}
for any integer $p\geq 1.$ Random measures satisfying the Ghirlanda-Guerra identities 
arise as the directing measures (or {\it determinators} in the terminology of \cite{SG2})
 of overlap matrices in the Sherrington-Kirkpatrick model  where they can be seen as 
 the asymptotic analogues of the Gibbs measure. The importance of the asymptotic point
 of view provided by these measures was brought to light in \cite{AA}, even though it was
 the Aizenman-Contucci stochastic stability \cite{AC} and not the Ghirlanda-Guerra identities that 
 played the main role there. However, subsequently, such random measures satisfying
 the Ghirlanda-Guerra identities played an equally important role in the results of \cite{PGG}
 and to a lesser extent of \cite{Tal-New}.
 
Our main result is based on a simple observation which extends the invariance theorem from \cite{PGG}.
Consider independent symmetric Bernoulli random variables $(\eps_\alpha)_{\alpha\in A}$ 
(taking values $\pm1$ with probability $1/2$) and for $t\in\Reals$ let us define a random measure 
$\mu_t=\sum_{\alpha\in  A}w_{\alpha,t} \delta_{\xi_\alpha}$ 
with weights defined by the random change of density
\begin{equation}
w_{\alpha,t}=\frac{w_\alpha \exp t\eps_\alpha}{\sum_{\gamma \in A} w_\gamma \exp t\eps_\gamma}.
\label{density}
\end{equation}
and as in (\ref{Av}) let us denote the average with respect to this measure by 
\begin{equation}
\la f \ra_t = \sum_{\alpha^1,\ldots, \alpha^n} w_{\alpha^1,t}\cdots w_{\alpha^n,t} 
f(\alpha^1,\ldots,\alpha^n).
\label{Avt}
\end{equation}
The following holds.
\smallskip

\begin{theorem}\label{Th1}
If a random measure $\mu$ satisfies the Ghirlanda-Guerra identities (\ref{GG}) then 
for any $t\in\Reals$, any $n\geq 2$ and any bounded function $f$ of the overlaps 
on $n$ configurations we have $\E\la f \ra_t = \E \la f \ra.$
\end{theorem}
\smallskip

The main difference here is that the result holds for all $t\in\Reals$ compared to
$|t|<1/2$ as stated in Theorem 4 in \cite{PGG} which was sufficient for the main
argument there.  However, letting $t$ go to infinity we now obtain the following 
new invariance property. Let $\eta_\alpha = (\eps_\alpha+1)/2$ be independent random 
variables, now taking values $1$ and $0$ with probability $1/2$ and let 
$\mu^\prime=\sum_{\alpha\in  A}w_{\alpha}^\prime \delta_{\xi_\alpha}$ 
be the random measure defined by the change of density
\begin{equation}
w_{\alpha}^\prime
=
\frac{w_\alpha \eta_\alpha}{\sum_{\gamma \in A} w_\gamma \eta_\gamma}.
\label{density2}
\end{equation}
In other words, we randomly delete half of the point in the support of measure $\mu$
and renormalize the weights to define a probability measure $\mu^\prime$ on the
remaining points. The denominator in (\ref{density2}) is non-zero with probability one 
since it is well-known that unless the measure $\mu$ is concentrated at $0\in B$ (a case
we do not consider)
it must have infinitely many different points in the support in order to satisfy (\ref{GG}).
Let us define by $\la f\ra^\prime$ the average with respect to $\mu^\prime$.
\smallskip

\begin{theorem} \label{Th2} (Deletion invariance)
If a random measure $\mu$ satisfies the Ghirlanda-Guerra identities (\ref{GG}) then 
$\E\la f \ra^\prime = \E \la f \ra.$
\end{theorem}
\smallskip

{\it Remark 1.}
In particular, this implies that the measure $\mu^\prime$ also satisfies the Ghirlanda-Guerra
identities (\ref{GG})  and, hence, we can repeat the random deletion procedure as many times
as we want. This means that the deletion invariance also holds with 
random variables $(\eta_\alpha)$ taking values $1$ and $0$ with probabilities 
$1/2^s$ and $1-1/2^s$ correspondingly, for any integer $s\geq 1$. 

{\it Remark 2.} It is well-known that invariance for the averages as in Theorem \ref{Th2} implies
that the random measures $\mu$ and $\mu^\prime$ have the same distribution, up to rotations.
Let $(w_\ell)_{\ell\geq 1}$ be the weights $(w_\alpha)$ arranged in the non-increasing order and
let $(\xi_\ell)$ be the points $(\xi_\alpha)$ rearranged accordingly, so that 
$\mu=\sum_{\ell\geq 1} w_\ell \delta_{\xi_\ell}$. Similarly, let
$\mu^\prime=\sum_{\ell\geq 1} w_\ell^\prime \delta_{\xi_\ell^\prime}.$ 
Then arguing as at the end of the proof of Theorem 4 in \cite{PGG}  (or Lemma 4 in \cite{PDS})
one can show that
\begin{equation}
\bigl((w_\ell)_{\ell\geq 1},(\xi_\ell\cdot \xi_{\ell'})_{\ell,\ell'\geq 1}\bigr)
\stackrel{d}{=}
\bigl((w_\ell^\prime)_{\ell\geq 1},(\xi_\ell^\prime\cdot \xi_{\ell'}^\prime)_{\ell,\ell'\geq 1}\bigr)
\end{equation}
which means that up to rotations the configurations of the random measures $\mu$ and $\mu'$ 
have the same distributions.
\smallskip

\noindent\textbf{Proof of Theorem \ref{Th1}.}
Suppose that $|f|\leq 1.$ Let $\varphi(t) = \E \la f\ra_t$ and by symmetry we only need to consider $t\geq 0$. 
Given configurations $\alpha^1,\alpha^2,\ldots$ let us denote
$$
D_n = \eps_{\alpha^1}+\ldots+\eps_{\alpha^n} - n \eps_{\alpha^{n+1}}
$$
and a straightforward computation shows that $\varphi'(t) = \E\la f D_n \ra_t$ and
similarly for all $k\geq 1,$
$$
\varphi^{(k)}(t) = \E\la f D_n \cdots D_{n+k-1}\ra_t.
$$
It was proved in Theorem 4 in \cite{PGG} (a more streamlined proof  was given in Theorem 6.3 
in \cite{Tal-New}) that if the measure $\mu$ satisfies 
the Ghirlanda-Guerra identities (\ref{GG}) then 
\begin{equation}
\varphi^{(k)}(0)=0 \,\,\mbox{ for all } \,\, k\geq 1. 
\label{derzero}
\end{equation}
It is also easy to see that $|D_n|\leq 2n$ so that for all $t,$
\begin{equation}
|\varphi^{(k)}(t)|\leq 2^k n(n+1)\cdots(n+k-1).
\label{derbound}
\end{equation}
This is all one needs to show that if 
\begin{equation}
\varphi(t) = \varphi(0) \,\,\mbox{ and }\,\,
\varphi^{(k)}(t)=0 \,\,\mbox{ for all }\,\, k\geq 1
\label{induction}
\end{equation}
holds for all $t\leq t_0$ for some $t_0\geq 0$ then it also holds for all $t<t_0+1/2.$ 
This will finish the proof of the theorem since by (\ref{derzero}) this holds for $t_0=0$.
Take any $k\geq 0.$ Using  (\ref{derbound}) and (\ref{induction}) for $t=t_0$ and using Taylor's expansion 
for a function $\varphi^{(k)}(t)$ around the point $t=t_0$ we get for any $m\geq 1$
\begin{eqnarray*}
&&
|\varphi^{(k)}(t)-\varphi^{(k)}(t_0)|
\leq
\sup_{t_0\leq s\leq t} \frac{|\varphi^{(k+m)}(s)|}{m!} |t-t_0|^m
\leq
\frac{ 2^{k+m} n(n+1)\cdots(n+k+m-1)}{m!} |t-t_0|^m.
\end{eqnarray*}
If $|t-t_0|<1/2$ then letting $m\to\infty$ proves that $\varphi^{(k)}(t)=\varphi^{(k)}(t_0)$ for all $k\geq 0$
and therefore (\ref{induction}) holds for all $t<t_0+1/2$. \qed

\noindent\textbf{Proof of Theorem \ref{Th2}.}
Let $I=\{\alpha\in A: \eps_\alpha=1\}$ and let 
$$
Z_t = \sum_{\alpha \in I} w_\alpha + e^{-2t} \sum_{\alpha\in I^c} w_\alpha
$$ 
so that
\begin{equation}
w_{\alpha,t} = \frac{w_\alpha}{Z_t} \bigl( I(\alpha\in I)+e^{-2t} I(\alpha\in I^c)\bigr).
\label{groups}
\end{equation}
Then the sum on the right hand side of  (\ref{Avt}) can by broken into $2^n$ groups
depending on which of the  indexes $\alpha^1,\ldots,\alpha^n$ belong to $I$ or its
complement $I^c,$ for example, the terms corresponding to all indices belonging to $I$ will give
\begin{equation}
\frac{1}{Z_t^n}\sum_{\alpha^1,\ldots, \alpha^n\in I} w_{\alpha^1}\cdots w_{\alpha^n} 
f(\alpha^1,\ldots,\alpha^n).
\label{Avt2}
\end{equation}
This sum is bounded by one and when $t\to+\infty$ it obviously converges to $\la f \ra^\prime$ while
the sums corresponding to other groups, when at least one of the indices belongs to $I^c$,
will converge to zero because of the factor $e^{-2t}$ in (\ref{groups}). By dominated convergence
theorem we get convergence of expectations. 
\qed

\end{document}